\documentclass[12pt]{article}
\begin{document}

\textwidth15.6cm
\textheight25.7cm
\normalbaselineskip=12pt
\normalbaselines
\parindent0.8cm
\hoffset-1cm
\voffset-3cm
\pagestyle{empty}
%%%%%%%%%%%%%%%%%%%%%%%%%%%%%%%%%%%%%%%%%%%%%%%%%%%%%%%%%%%%%%%%%
%
%  Here the style of numbering and referencing is given.
%  Do not change!!!!!
%%%%%%%%%%%%%%%%%%%%%%%%%%%%%%%%%%%%%%%%%%%%%%%%%%%%%%%%%%%%%%
\catcode `\@=11
\@addtoreset{equation}{section}

\def\theequation{\arabic{section}.\arabic{equation}}
          % if you want equations to be numbered by section

\def\section{\@startsection {section}{1}{\z@}{-3.5ex plus -1ex minus
     -.2ex}{2.3ex plus .2ex}{\normalsize\bf}}
\def\subsection{\@startsection{subsection}{2}{\z@}{-3.25ex plus -1ex
minus
 -.2ex}{1.5ex plus .2ex}{\normalsize\bf}}
          % correct font size for section/subsection titles

\def\thebibliography#1{\section*{References\markboth
  {REFERENCES}{REFERENCES}}\list
  {[\arabic{enumi}]}{\settowidth\labelwidth{[#1]}\leftmargin\labelwidth
  \advance\leftmargin\labelsep
  \usecounter{enumi}}
  \def\newblock{\hskip .11em plus .33em minus -.07em}
  \sloppy
  \sfcode`\.=1000\relax}
 \let\endthebibliography=\endlist
% numbering of references as ``[3] Author''

\catcode `\@=12

%%%%%%%%%%%%%%%%%%%%%%%%%%%<Personal macros>%%%%%%%%%%%%%%%%%%%%%%%%
%%%%%%%%%%%
%%%%%%%%%%%%%%%%%%%%%%%%%%%%%%%%%%%%%%%%%%%
\def\g{\mbox{\bf g\,}}
\def\R{\mbox{$\cal R$\,}}
\def\F{\mbox{$\cal F$}}
\def\1{{\bf 1}}
\def\ot{\otimes}
\def\id{\mbox{id}}
\def\A{\mbox{$\cal A$}}
\def\P{\mbox{$\cal P_{\mp}{}$}}
%Macro for the black-board bold
\def\b#1{{\mathbb #1}}
\newcommand{\tr}{\triangleright_h}
\newcommand{\trc}{\triangleright}
\newtheorem{prop}{Proposition}
\newtheorem{lemma}{Lemma}
\newtheorem{theorem}{Theorem}
\newtheorem{corollary}{Corollary}

%%%%%%%%%%%%%%%%%%%%%%%%%%%%%%%%%%
%
%   HERE  IS  THE START
%
%%%%%%%%%%%%%%%%%%%%%%%%%%%%%%%%%%%%%%%%%%%%%%%%

\vspace*{2.5cm}
\noindent
{ \bf EMBEDDING $q$-DEFORMED HEISENBERG \\ 
ALGEBRAS INTO UNDEFORMED ONES}
\footnote{Talk given at the XVI-th Workshop On Geometric Methods 
in Physics, July 1997, Bialowieza (Poland). To appear
in the Proceedings, to be published on Rep. Math. Phys.}
\vspace{1.3cm}\\
\noindent
\hspace*{1in}
\begin{minipage}{13cm}
Gaetano Fiore  \vspace{0.3cm}\\
 Dip. di Matematica e Applicazioni, Universit\`a di Napoli\\
 \hspace*{4cm} and \hfill\\
 Istituto Nazionale di Fisica Nucleare, Sezione di Napoli\\
 Italy 
\end{minipage}

\vspace*{0.5cm}

\begin{abstract}
\noindent
Any deformation of a Weyl or Clifford algebra can  be realized 
through some change of generators in the undeformed algebra.
Here we briefly describe and motivate 
our systematic procedure for constructing
all such changes of generators for those particular
deformations where the original algebra is covariant under
some Lie group and the deformed algebra is covariant
under the corresponding quantum group.
\end{abstract}

% section 1
\section{\hspace{-4mm}.\hspace{2mm} INTRODUCTION}

Weyl and Clifford algebras are at the hearth of quantum physics.
One may ask  if deforming them, i.e. deforming their 
defining commutation relations, yields new physics 
\cite{flato}, or at least may be useful to better describe
some peculiar systems in conventional quantum physics.
This question can be divided into an algebraic and a 
representation-theoretic
subquestions. Roughly speaking, the first is: 
is there a formal realization of the elements of the deformed
algebra in terms of elements of the undeformed algebra? 
The answer is affirmative \cite{gerst,ducloux,pillin}
but in general the realization is not explicitly known. 
The second 
subquestion is: do also the corresponding representation theories
coincide? One can already see in some simple model that the answer is
negative, but in the general case, up to our knowledge, 
the relation between the two is an open question.

We introduce the notions of a deformed algebra and of a
deforming map first on a simplest toy model, the 1-dim Weyl algebra
$\A$.
$\A$ is generated by $\1, a,a^+$ fulfilling
\begin{equation}
a \,a^+=\1+a^+a  \qquad\qquad \1 b=b \1=b ,
\end{equation}
$b\in\A$. As a {\it deformation}
$\A_h$ of $\A$ ($h$ is the `deformation parameter') we consider the
algebra
generated by $\1_h,\tilde A,\tilde A^+$ fulfilling the relations
\begin{equation}
\tilde A\,\tilde A^+=\1_h+ e^h\tilde A^+\,\tilde A
\qquad\qquad \1_h B=B\1_h=B,                                
\label{proto}
\end{equation}
$B\in\A_h$;
when $h\rightarrow 0$ the second relations go to the first
if we identify in the limit $\1_h,\tilde A,\tilde A^+$ with
$\1,a,a^+$.

Can we realize $\1_h,\tilde A,\tilde A^+$ within $\A[[h]]$ (the ring of
formal power series in the unknown $h$ and with coefficients in $\A$),
in other words as `functions' 
of $h,a,a^+$ reducing to $\1,a,a^+$ in the limit? Yes. Let  $n:=
a^+\,a$,
$q:= e^h$, $(x)_q:= \frac{q^x-1}{q-1}$; if we define \cite{zachos}
\begin{equation}
A:=a \sqrt{\frac{(n )_q}{n }}
\qquad\qquad
A^+:=\sqrt{\frac{(n )_q}{n }}a^+,   \label{traproto}
\end{equation}
it is easy to show that $\1, A,A^+$ indeed fulfil the `deformed 
commutation
relations' (DCR)  (\ref{proto}): in other words
$\1,A,A^+$ realize $\1_h,\tilde A,\tilde A^+$. At lowest order in $h$
one 
finds
$A=a+O(h)$, $A^+=a^++O(h)$, as required.
By definition, a {\it deforming
map} $f$ is an algebra isomorphism $f:\A_h\rightarrow \A[[h]]$ over
${\bf C}[[h]]$ reducing to the identity in the limit $h=0$.
We can obtain one by setting $f(\tilde A):=A$,
$f(\tilde A^+):=A^+$ and extending its action on the whole $\A_h$
imposing $f(\alpha\beta)=f(\alpha)f(\beta)$ and linearity.

Here we shall deal with a particular class of deformations of 
multidimensional
Weyl algebras or Clifford algebras (their fermionic counterparts).
The undeformed algebra is covariant under some Lie algebra $\g$
and the deformed one under the quantum group \cite{dr2} $U_h\g$.
The undeformed algebra $\A$ is generated by $\1, a^i,a^+_j$ fulfilling
\begin{eqnarray}
&&[a^i\, , \,a^j]_{\mp}            = 0 \nonumber\\~
&&[a^+_i\, , \, a^+_j]_{\mp}  = 0 \label{ccr}\\~
&&[a^i\, , \, a^+_j]_{\mp}       = \delta_j^i\1
\nonumber
\end{eqnarray}
(the $\mp$ sign denotes commutators and
anticommutators and refers to Weyl and Clifford algebras respectively)
and transforms under the action $\trc$ of $\g$ according to some law
\begin{equation}
    x\trc a^+_i=\rho(x)^j_ia^+_j\qquad\qquad
      x\trc a^i=\rho(Sx)_j^ia^j;
\end{equation}
here $x\in\g$, $Sx=-x$ and $\rho$ denotes
some matrix representation of \g. Clearly $a^i$ belong to a
represenation 
of $\g$ which is the contragradient of the $a^+_i$ one. The action
$\trc$ is extended to products of the
generators using the standard rules of tensor product representations
(technically speaking, using the coproduct $\Delta$ of 
the universal enveloping algebra $U\g$), and then
linearly to all of $\A$; this is possible because the action
of $\g$ is manifestly compatible with the commutation relations 
(\ref{ccr}).
 The same formulae, where $S$ now denotes
the antipode of $U\g$,  give also the standard extension 
of $\trc$ to $x\in U\g$.
 
The corresponding
deformed algebra $\A_h$ is generated by $\1_h,\tilde A^+_i,\tilde A^i$
fulfilling DCR  which, in the simplest case of $\rho$ being
the defining fundamental representation of $\g$, take the form
\cite{puwo,wezu,cawa}
\begin{eqnarray}
&& \P^{ij}_{hk}  \tilde A^k \tilde A^h      = 0 \nonumber\\~
&& \P_{ij}^{hk}  \tilde A^+_h \tilde A^+_k  =0 \label{dcr}\\~
&& \tilde A^i\tilde A^+_j= \delta^i_j\,\1 \pm q \hat R^{ih}_{jk}
\tilde A^+_h \tilde A^k;
\nonumber
\end{eqnarray}
$\A_h$ transforms under the action $\tr$ of $U_h\g$ according to the law
\footnote{These $\A_h$ should not be confused with the
celebrated Biedenharn-Macfarlane-Hayashi-Kulish $q$-oscillator 
(super)algebrae \cite{mac}, whose generators 
$\alpha^i,\alpha^+_j$ fulfil ordinary (anti)commutation relations,
except for the $q$-(anti)commutation relations 
$\alpha^i\alpha^+_i\mp q^2\,\alpha^+_i\alpha^i=1$, and are {\it not} 
$U_h\g$-covariant (in spite of the fact that they are usually used
to construct a generalized Jordan-Schwinger realization of $U_h\g$).}
\begin{equation}
x\,\tr\tilde A^+_i = \rho_h{}^j_i(x)A^+_j
\qquad\qquad
x\,\tr\tilde A^i = \rho_h{}^i_j(S_h x)A^j.
\label{qtrans}
\end{equation}
Here $x\in U_h\g$, $S_h$ is the antipode of $U_h\g$,
$\rho_h{}$ the quantum group deformation of $\rho$,
$\hat R$ the braid matrix of $U_h\g$ in the representation $\rho_h{}$,
and finally ${\cal P}_-,{\cal P}_+$ are the 
$U_h\g$-covariant deformations of the
antisymmetric and symmetric projectors, in the form of
polynomials in $\hat R$; for instance,  when $\g=sl(N)$
$\P= (q+q^{-1})^{-1} ( q^{\pm 1}\1\mp\hat R)$. The upper and lower sign
refer to Weyl and Clifford algebras respectively.
$\tilde A^i$ belong to a representation of $U_h\g$ which is the 
quantum group contragradient of the $\tilde A^+_i$ one. The action
$\tr$ is extended to products of the
generators using the coproduct $\Delta_h$ of $U_h\g$, see eq.
(\ref{modalg}) below, and then
linearly to all of $\A_h$; this is possible because the action
of $U_h\g$ is compatible with the commutation relations (\ref{dcr}).

Is there a realization of $\1_h,\tilde A^i,\tilde A^+_j$ within
$\A[[h]]$,
or, in other words, a deforming map $f:\A_h\rightarrow \A[[h]]$? Yes. 
The affirmative answer is based on the vanishing of the second
Hochschild 
cohomology
group of any Weyl algebra
\cite{gerst,ducloux,pillin}; this allows to prove the existence of $f$ 
without
however providing an explicit construction. 
The argument is valid not only for deformations of the
type (\ref{dcr}), but for {\it any} kind of deformation of $\A$.

In \cite{fiormp} we have suggested a systematic and explicit 
constructing procedure of deforming maps for the class of Weyl 
and Clifford algebras described above; the procedure is
based on $U_h\g$-covariance and the socalled
Drinfel'd twist \cite{dr3}. In the sequel we briefly describe it.
We shall motivate our physical interest in these deforming maps in the
last section.

% section 2
\section{\hspace{-4mm}.\hspace{2mm} THE CONSTRUCTING PROCEDURE }

If $\alpha\in\A[[h]]$ is any
element of the form $\alpha=\1+O(h)$ and $f$ is a deforming map,
one can obtain a new one $f_{\alpha}$  by the inner automorhism
\begin{equation}
f_{\alpha}(\cdot):=\alpha f(\cdot)\alpha^{-1};
\label{inner}
\end{equation}
 actually the vanishing of the 
first Hochschild cohomology group of $\A$ implies that
{\it all} deforming maps can be obtained from one
in this manner. Therefore our problem is reduced to 
finding a particular one, what we are going to describe below.

The other essential ingredients of our construction procedure are:
\begin{enumerate}

\item \g, a simple Lie algebra.
\item The cocommutative Hopf algebra $H\equiv(U\g,\cdot,\Delta,
      \varepsilon,S)$ associated to 
      $U\g$; $\cdot,\Delta,\varepsilon,S$ denote the product,
      coproduct, counit, antipode.
      We shall use the Sweedler's notation 
      $\Delta(x)\equiv x_{(1)}\ot x_{(2)}$: at the rhs a sum 
      $\sum_i x^i_{(1)}\ot x^i_{(2)}$ of many terms is implicitly 
      understood.
\item The quantum group \cite{dr2} $H_h\equiv(U_h\g,\bullet,\Delta_h,
      \varepsilon_h,S_h,\R)$. $\bullet,\Delta_h,\varepsilon_h,S_h$
denote
      the deformed product, coproduct, counit, antipode, $\R$ the
      quasitriangular structure. We shall use the Sweedler's notation 
      (with barred indices)
      $\Delta_h(x)\equiv x_{(\bar 1)}\ot x_{(\bar 2)}$.
\item An algebra isomorphism\cite{dr3} 
      $\varphi_h:U_h\g\rightarrow U\g[[h]]$ over ${\bf C}[[h]]$, namely
      $\varphi_h(x\bullet y)=\varphi_h(x)\cdot\varphi_h(y)$.
\item A corresponding Drinfel'd twist\cite{dr3} 
      $\F\equiv\F^{(1)}\!\ot\!\F^{(2)}\!=\!\1^{\ot^2}\!\!
      +\!O(h)\in U\g\![[h]]^{\ot^2}$:
      \begin{equation}
      (\varepsilon\ot \id)\F=\1=(\id\ot \varepsilon)\F,
      \qquad\: \:\Delta_h(a)=(\varphi_h^{-1}\ot \varphi_h^{-1})\big
      \{\F\Delta[\varphi_h(a)]\F^{-1}\big\};
      \end{equation}
      the last formula means that, up to the isomorphism $\varphi_h$,
      $\Delta_h$ is
      related to $\Delta$ by a similarity transformation.
\item $\gamma':=\F^{(2)}\cdot S\F^{(1)}$ and 
      $\gamma:=S\F^{-1(1)}\cdot \F^{-1(2)}$. Up to the isomorphism 
      $\varphi_h$,
      $S_h$ and its inverse are related to $S$ by similarity 
      transformations
      involving resp. $\gamma$ and $\gamma'$.
\item The Jordan-Schwinger algebra homomorphism 
      $\sigma:U\g[[h]]\rightarrow\A[[h]]$, defined on the 
      generators by 
      \begin{equation}
      \sigma(\1_{U\mbox{\footnotesize \bf g}})=\1\qquad\qquad
      \sigma(x):=
      \rho(x)^i_ja^+_ia^j
      \end{equation}
      $x\in\g$, and extended to the whole $U\g[[h]]$ as an algebra
      homomorphism, $\sigma(yz)=\sigma(y)\sigma(z)$
      and $\sigma(y+z)=\sigma(y)+\sigma(z)$. This is consistent
      because $\sigma([x,y])=[\sigma(x),\sigma(y)]$. In the $su(2)$
      $\sigma$ takes the well-known form
      \begin{equation}
      \sigma(j_+)=a^+_{\uparrow}a^{\downarrow},\qquad\qquad
      \sigma(j_-)=a^+_{\downarrow}a^{\uparrow},\qquad\qquad
      \sigma(j_0)=\frac 12(a^+_{\uparrow}a^{\uparrow}-
      a^+_{\downarrow}a^{\downarrow}).
      \label{homo}
      \end{equation}
\item The deformed Jordan-Schwinger algebra homomorphism 
      $\sigma_h:U_h\g\rightarrow\A[[h]]$, defined by
      $\sigma_h:=\sigma\circ\varphi_h$.
\item The $*$-structures $*,*_h,\star,\star_h$ in $H,H_h,\A,\A_h$, if 
      $\A,\A_h$ are $*$-algebras transforming respectively under the
Hopf
      $*$-algebras $H,H_h$ with the compatibility condition
      \begin{equation}
      (x\,\tr a)^{\star_h}=S_h^{-1}(x^{*_h})\tr a^{\star_h}.
      \label{condstar}
      \end{equation}
\end{enumerate}

\subsection{\hspace{-5mm}.\hspace{2mm} Constructing the Quantum Group 
Action and the Generators $A^i,A^+_j$}

Since we know that a deforming map exists, although we cannot write it 
explicitly we can say that it must be possible to realize $\tr$ on 
$\A[[h]]$,
instead of $\A_h$. Our first step is to guess such a realization. 
This requires fulfilling 
\begin{equation}
(xy)\tr a=x\tr(y\tr a) \qquad\qquad\qquad
x\tr(ab)=(x_{(\bar 1)}\tr a)(x_{(\bar 2)}\tr b)
\label{modalg}
\end{equation}
for any $x,y\in U_h\g$, $a,b\in \A_h$; these are the conditions 
characterizing
a module algebra. 
There is a simple way to find such a realization, namely by setting
\begin{equation}
x\tr a := \sigma_h(x_{(\bar 1)}) a 
\sigma_h(S_h x_{(\bar 2)});
\label{defprop}
\end{equation}
it is easy to check that (\ref{modalg}) are indeed fulfilled using the 
basic
axioms characterizing the coproduct, counit, antipode in a generic
Hopf algebra.
The guess has been suggested by the cocommutative case, where the
same conditions  and realization are obtained
for $U\g,\A,\trc$ if in the two previous formulae we just erase the 
suffix ${}_h$ and replace $\Delta_h(x)\equiv x_{(\bar 1)}\ot x_{(\bar
2)}$
with the cocommutative coproduct $\Delta(x)\equiv x_{(1)}\ot x_{(2)}$.
 
Our second step is to realize elements $A^i,A^+_j\in\A[[h]]$ that
transform under (\ref{defprop}) as $\tilde A^i,\tilde A^+_j$
in (\ref{qtrans}). Note that $a^i,a^+_j$ do {\it not} transform
in this way. In Ref. \cite{fiojpa} we proved that the following objects
do:
\begin{equation}
\begin{array}{lll}
A_i^+ &:= & u\, \sigma(\F^{(1)})a_i^+
\sigma(S \F^{(2)}\gamma)\, u^{-1} \cr
A^i&:= & v\, \sigma(\gamma'S \F^{-1(2)})a^i
\sigma(\F^{-1(1)}) v^{-1};          \cr                
\end{array}
\label{def3}
\end{equation}
the result holds for any choice of $g$-invariant elements $u,v=\1+O(h)$ 
in $\A[[h]]$, in particular for $u=v=\1$.

The third step is to fix $u,v$ in such a way that the DCR are fulfilled.
One can easily show that the DCR may fix at most the product $u v^{-1}$.
For the explicit case considered in (\ref{dcr}) we proved  in Ref. 
\cite{fiormp} that the
DCR are indeed fulfilled by taking $u v^{-1}$ 
\begin{eqnarray}
u v^{-1} &\! = \!&\frac{\Gamma(n+1)}{\Gamma_{q^2}(n+1)}          \\
u v^{-1} &\! =\! &\left(\frac{1\!+\!q^{N-2}}2\right)^{-n}
\frac{\Gamma\left[\frac 12\left(n\!+\!\frac{N}2\!+\!1\!-\!l
\right)\right]\,\Gamma\left[\frac 12\left(n\!+\!
\frac{N}2\!+\!1\!+\!l\right)\right]}{\Gamma_{q^2}
\left[\frac 12\left(n\!+\!\frac{N}2\!+\!1\!-\!l
\right)\right]\,\Gamma_{q^2}\left[\frac 12
\left(n\!+\!\frac{N}2\!+\!1\!+\!l\right)\right]},
\end{eqnarray}
respectively if $\g=sl(N),so(N)$. 
Here $\Gamma$ is Euler's $\gamma$-function, $\Gamma_{q^2}$  its 
$q$-deformation characterized by $\Gamma_{q^2}(x+1)=(x)_{q^2}
\Gamma_{q^2}(x)$,
$n:= a^ia^+_i$, and in the $\g=so(N)$ case
we have enlarged for convenience $\A[[h]]$ by the introduction of the 
square root $l:=\sqrt{\sigma({\cal C})}$,  ${\cal C}$ being the
quadratic
casimir of $Uso(N)$. We stress that the above solutions regard the case 
of $\rho$ being  the defining representation of $\g$.
We have yet no  formula yielding the right 
$u v^{-1}$, if any, necessary to fulfil the DCR in the general case.
However it is important to note that in general the DCR translate
into  conditions on $u v^{-1}$ 
where  the Drinfel'd twist $\F$ appears only through
the socalled `coassociator'
\begin{equation}
\phi:=[(\Delta\ot \id)(\F^{-1})](\F^{-1}\ot \1)(\1 \ot \F)
[(\id \ot \Delta)(\F).
\label{defphi}
\end{equation}
$\phi$ is explicitly known, unlike $\F$, for which up to now
there is an existence proof but no explicit expression.
This makes the above conditions explicit and allows to search the
explicit form of $u v^{-1}$ in the general case, if it exists.

Finally, the residual freedom in the choice of $u,v$ is partially fixed
if $H,H_h,\A,\A_h$ are matched (Hopf) $*$-algebras
and make the additional requirement that $\star$ realizes in $\A[[h]]$
the $\star_h$ of $\A_h$. For instance, if $(a^i)^{\star}=a^+_i$ and
$h\in {\bf R}$ this means
\begin{equation}
(A^i)^{\star}=A^+_i,
\end{equation}
and is fulfilled if we take $u=v^{-1}$.

In the $\g=sl(2)$ case, with $\rho$ being the fundamental
representation,
the knowledge of $(\rho\otimes\id)\F$ is sufficient to determine the 
$A^i,A^+_i$ of
formulae (\ref{def3}) completely. Denoting $n^i=a^ia^+_i$
(with {\it no} sum over $i$),  $i=\uparrow, \downarrow$
and taking $u=v^{-1}$ one finds for the deformed Weyl algebra
\begin{equation}
\begin{array}{rclcrcl}
A^+_{\uparrow} & = &\sqrt{(n^{\uparrow})_{q^2}\over n^{\uparrow}}
q^{n^{\downarrow}}a^+_{\uparrow} &\qquad
\qquad A^+_{\downarrow} & = &
\sqrt{(n^{\downarrow})_{q^2}\over n^{\downarrow}}
a^+_{\downarrow}  \nonumber \\
A^{\uparrow} & = & a^{\uparrow}\sqrt{(n^{\uparrow})_{q^2}\over 
n^{\uparrow}}
q^{n^{\downarrow}} &\qquad
\qquad A^{\downarrow} & = &a^{\downarrow} 
\sqrt{(n^{\downarrow})_{q^2}\over n^{\downarrow}},
\end{array}
\label{lastb}
\end{equation}
and for the deformed Clifford one
\begin{equation}
\begin{array}{rclcrcl}
A^+_{\uparrow} & = & q^{-n^{\downarrow}}a^+_{\uparrow} &\qquad
\qquad A^+_{\downarrow} & = &
a^+_{\downarrow}  \nonumber \\
A^{\uparrow} & = &a^{\uparrow}q^{-n^{\downarrow}} &\qquad
\qquad A^{\downarrow} & = & a^{\downarrow}.
\end{array}
\label{lastf}
\end{equation}

In the case that the Hopf algebra $H_h$ is not a genuine quantum group,
but
a triangular one, the whole discussion simplifies in that one can take
trivial $u,v$, see \cite{fiojmp}.

Above we have determined in $\A[[h]]$ one particular realization 
$A^i,A^+_j$  and $\tr$
of the generators $\tilde A^i,\tilde A^+_j$ and of the quantum group 
action.
Its main feature is that 
the $\g$-invariant ground state $|0\rangle$ as well as the first excited
states $a^+_i|0\rangle$ of the classical Fock space representation
are also respectively $U_h\g$-invariant
ground state $|0_h\rangle$ and first excited
states $A^+_i|0_h\rangle$ of the deformed Fock space representation.

According to eq. (\ref{inner}) all the other realizations are of the
form
\begin{equation}
A^{\alpha\,i}=\alpha\,A^i\alpha^{-1}\qquad\qquad
A^{+}_{\alpha\,i}= \alpha\, A^+_i\, \alpha^{-1},
\label{aa'}
\end{equation}
with $\alpha=\1+O(h)\in\A[[h]]$.
They are manifestly covariant under the realization 
$\trc_{h,\alpha}$ of the $U_h\g$-action defined by
\begin{equation}
x\,\trc_{h,\alpha}\,a\: :=\:\alpha\;\sigma_h(x_{(\bar 1)}) 
\, a\,\sigma_h(x_{(\bar 2)})\;\alpha^{-1}.
\end{equation}
For these realizations the deformed ground state in the Fock space 
representation
reads $|0_h\rangle=\alpha|0\rangle$; if $\alpha|0\rangle\neq |0\rangle$ 
the $\g$-invariant ground state and first excited
states of the classical Fock space representation
do not coincide with their deformed counterparts.

% section 3
\section{\hspace{-4mm}.\hspace{2mm} CLASSICAL VS. QUANTUM INVARIANTS}
\label{inva}
We have introduced two actions on $\A[[h]]$: 
\begin{equation}
\trc: U\g \times \A[[h]]\rightarrow \A[[h]], \qquad\qquad
\tr: U_h\g \times \A[[h]]\rightarrow \A[[h]].
\end{equation}
Their respective invariant subalgebras
$\A^{inv}[[h]],\A^{inv}_h[[h]]$ are defined by
\begin{equation}
\A^{inv}_h [[h]] : =\{I\in \A[[h]]\:\: |\:\: x\tr I=\varepsilon_h(x) 
I\qquad\forall x\in U_h\g\}
\label{def7}
\end{equation}
and by the analogous equation where all suffices ${}_h$ are erased.
What is the relation between them? It is easy to prove that
\cite{fiormp}
\begin{equation}
\A^{inv}_h[[h]]=\A^{inv}[[h]].
\label{prop}
\end{equation}
In other words invariants under the $\g$-action $\trc$ are also
$U_h\g$-invariants under $\tr$, and conversely, although in general
$\g$-covariant objects (tensors) and $U_h\g$-covariant ones do not
coincide in general!

Let us introduce 
in the vector space $\A^{inv}[[h]]=\A^{inv}_h[[h]]$ bases $I^1,I^2,...$ 
and $I^1_h,I^2_h,...$ consisting
of polynomials respectively in $a^i, a^+_j$ and $A^i,A^+_j$. It is 
immediate
to realize that we can choose the polynomials homogeneous, since
$\trc,\tr$ act linearly without changing their degrees.
Explicitly,
\begin{equation}
\begin{array}{lll}
I^1=a^+_ia^i                &\qquad \qquad &    I^1_h=A^+_iA^i  \\
I^2=d^{ijk}a^+_ia^+_ja^+_k  &\qquad \qquad &     
I^2_h=D^{ijk}A^+_iA^+_jA^+_k\\
I^3=d'_{kji}a^ia^ja^k       &\qquad \qquad &     
I^3_h=D'_{kji}A^iA^jA^k  \\
I^4=....                    &\qquad \qquad &    I^4_h=....   
\end{array}
\end{equation}
where the numerical coefficients  $d,d',...$ 
form $\g$-isotropic tensors and the numerical coefficients
$D,D'$ the corresponding $U_h\g$-isotropic tensors.
It is easy to show that $I^1_h\neq I^1$. In general 
$I^n_h\neq I^n$, although $I^n_h=I^n+O(h)$.
The propostion (\ref{prop}) implies in particular
\begin{equation}
I^n_h=g^n(\{I^m\},h)=k^n(\{a^i,a^+_j\},h).
\end{equation}
What do the `functions' $g^n,k^n$, i.e. the formal power series
in $h$ with coefficients respectively in $\A^{inv}$ and $\A$,
look like?

In Ref. \cite{fiormp} we have found universal formulae yielding
the $k^n$'s. The latter turn out to be highly non-polynomial functions, 
or more precisely in their power expansions in $h$
the degree in $a^i, a^+_j$ of the corresponding polynomial coefficient 
grows without bound with the power.
It is remarkable that in these universal formulae the twist $\F$
appears only through the coassociator $\phi$; therefore all the
$k^n$ can be worked out explicitly.

In the case that the Hopf algebra $H_h$ is not a genuine quantum group, 
but 
triangular, the coassociator as well as $u,v$ are trivial and one finds
$I^n_h=I^n$.

% section 3
\section{\hspace{-4mm}.\hspace{2mm} FINAL REMARKS, MOTIVATIONS 
AND CONCLUSIONS}

We have shown how one can can realize a deformed $U_h\g$-covariant
Weyl or Clifford algebra $\A_h$
within the undeformed one $\A[[h]]$. Given a representation $(\pi,V)$
of $\A$ on a vector space $V$, does it provide also a representation
of $\A_h$? In other words, can one interpret the elements of $\A_h$ 
as operators acting on $V$, if the elements of $\A$ are? If so,
which specific role play the elements $A^i,A^+_i$ of $\A[[h]]$?

In view of the specific example we have examined in ref. \cite{fiojmp}
the answer to the first question seems to be always positive, whereas
the
converse statement is wrong: there are more (inequivalent)
representations of the deformed algebra than representations of
the undeformed ones. This may seem surprising, but is not really a
paradox,
since the limit $h\rightarrow 0$ is smooth for the deforming maps
and their inverses only in a $h$-formal-power sense. 
Of course, we are especially interested in Hilbert
space representations of $*$-algebras.
In Ref. \cite{fiojmp} we checked that in the operator-norm topology
$f^{-1}$ is ill-defined on all but one `deformed' representation.
Roughly speaking, the reason is that the `particle-
number' observables $n^i$, which enter the transformation $f$ (see
e.g. (\ref{lastb})) are unbounded operators, therefore even for very 
small $h$ the effect of the transformation on their large-eigenvalue 
eigenvectors can be so large to `push' the latter out of the domain of
definition of the operators in $f^{-1}(\A)$.

We are especially interested in the case of $*$-algebras admitting Fock 
space representations. The results presented in the previous paragraphs
could in principle be applied to models in quantum field
theory or condensed matter physics by choosing representations
$\rho$ which are the direct sum of many copies of the same
fundamental representation $\rho_d$; this is what we have 
addressed in Ref. \cite{fiojpa}. The different copies would 
correspond respectively to different space(time)-points or
crystal sites. 

One important issue is if $U_h\g$-covariance necessarily
implies exotic particle statistics. In view of what we have said
the answer is no. At least for compact $\g$ and $U_h\g$ ($h$ is real),
the undeformed Fock space representation, which allows a
`Bosons \& Fermions' particle interpretation, carries also a
representation of the deformed one. Next point is the role
of the operators $A^i,A^+_j$. Quadratic commutation relations
of the type (\ref{dcr}) mean  that $A^+_i,A^i$ act as creators and
annihilators of some excitations; a glance at (\ref{def3}), (\ref{aa'})
shows that these are not the undeformed excitations,
but some `composite' ones. The last point is:
what could the latter be good for. As an Hamiltonian $H$ of the
system we can choose a simple combination of the $U_h\g$-invariants
$I^n_h$ of section \ref{inva}; the Hamiltonian is $U_h\g$-invariant and
has
a simple polynomial structure in the composite operators
$A^i,A^+_j$. $H$ is also $\g$-invariant, but has a highly non-polynomial
structure in the undeformed generators $a^i,a^+_j$ (it would be
tempting to understand what kind of physics it could describe!).
This suggests that the use of the $A^i,A^+_j$ instead of the $a^i,a^+_j$
should simplify the resolution of the corresponding dynamics.

%%%%%%%%%%% Replace the text between here and the next line begining
with
%%%%%%%%%%% by the text of your contribution.

%%%%%%%%%%% Insert your bibliography below %%%%%%%%%%


\begin{thebibliography}{99}

\footnotesize
\bibitem{flato} F.\ Bayen, M.\ Flato, C.\ Fronsdal, A.\ Lichnerowicz
and D.\ Sternheimer, {\em Ann. Phys.}, 111, 61, (1978).

\bibitem{gerst} M.\ Gerstenhaber,
%{\it On the Deformation of  Rings and Algebrae},
{\em  Ann.\ Math.} 79, 59, (1964).

\bibitem{ducloux} F.\ du Cloux, {\em Asterisque (Soc. Math. France)}
124-125, 129, (1985).

\bibitem{pillin} M.\ Pillin, 
%{\it On the Deformability of Heisenberg Algebrae},
{\em Commun. Math. Phys.} 180, 23, (1996). 

\bibitem{zachos} See e.g.: S.\ Vokos, C.\ Zachos, {\em Mod. Phys.
Lett.},
 A9, 1, (1994), and references therein.

\bibitem{dr2}
V.\ G.\ Drinfeld, 
Quantum groups, page 798 {\em in}
``Proceedings of the International
Congress of Mathematicians, Berkeley 1986'', Gleason Ed., Providence, 
(1987). M.\ Jimbo, {\em Lett. Math. Phys.} 10, 63, (1985). 

\bibitem{mac} L.\ C.\ Biedenharn,
{\em J. Phys. } A22, L873, (1989).
A.\ J.\ Macfarlane, 
{\em J.Phys.}  A22, 4581, (1989). T.\ Hayashi, 
{\em Commun. Math. Phys.} 127, 129, (1990).
 M.\ Chaichian and P.\ Kulish, {\em Phys. Lett.} B234, 72, (1990).


\bibitem{fiormp} G.\ Fiore, 
 Drinfel'd Twist and $q$-Deforming Maps for Lie Group
Covariant Heisenberg Algebras, e-print
q-alg/9708017, to appear in {\em Rev. Math. Phys}. 

\bibitem{fiojmp} G.\ Fiore, {\em J. Math. Phys.}  39, 3437, (1998).

\bibitem{dr3}
V.\ G.\ Drinfeld, 
%{\it Quasi Hopf Algebrae}, 
{\em Leningrad Math.\ J.} 1, 1419, (1990).

\bibitem{puwo} W. \ Pusz, S. \ L. \ Woronowicz, 
%{\it Twisted Second Quantization},
{\em Reports on Mathematical Physics} 27, 231, (1989).

\bibitem{wezu}
J.\ Wess and B.\ Zumino, 
%{\it Covariant Differential Calculus on the Quantum Hyperplane}
{\em Nucl. Phys. Proc. Suppl.} 18B, 302, (1991).

\bibitem{cawa}
U.\ Carow-Watamura, M.\ Schlieker and S.\ 
{\em Watamura, Z. Phys. C Part. Fields} 49,  439, (1991).


\bibitem{fiojpa} G.\ Fiore,
{\em J. Phys. A} 31, 5289, (1998).



%\bibitem{jqe} J. Q. Editor, ``Book Title, with the Initial Letter of
Each
%Major  Word Capitalized'', Publisher, City (1967)

%\bibitem{aa} A. Author, Article title, with only the first word having
an
%initial capital, {\em in}: ``Book Title'', A. Jones, ed., Publisher,
City
%(1967)

%\bibitem{bb} B. Bingo, {\em J. Irreprod. Results}  22, 123,
%(1988)

\end{thebibliography}
\end{document}